\documentclass{ifacconf}

\usepackage{amsmath,amssymb,amsfonts}
\usepackage{algorithmic}
\usepackage{graphicx}
\usepackage{textcomp}
\usepackage{xcolor}
\usepackage{epstopdf}   
\usepackage{natbib}        
\def\caly{{\mathcal Y}}

\def\calr{{\mathcal R}}

\def\calp{{\mathcal P}}

\def\cala{{\mathcal A}}

\def\calz{{\mathcal Z}}

\def\L2e{{\cal L}_{2e}}

\def\rea{\mathbb{R}}

\def\adj{\mbox{adj}}

\newtheorem{assumption}{Assumption}

\begin{document}
\begin{frontmatter}

\title{Adaptive observer for a nonlinear system with partially unknown state matrix and  delayed measurements\thanksref{footnoteinfo}} 

\thanks[footnoteinfo]{This work is supported by the Russian Science Foundation under grant 22-21-00499, https://rscf.ru/project/22-21-00499/}

\author[First]{Olga Kozachek} 
\author[Second]{Alexey Bobtsov}
\author[Third]{Nikolaev Nikolaev} 

\address[First]{ITMO University, Kronverksky av., 49, 197101, Saint Petersburg, Russia (e-mail: oakozachek@mail.ru)}
\address[Second]{ITMO University, Kronverksky av., 49, 197101, Saint Petersburg, Russia (e-mail:  bobtsov@mail.ru)}
\address[Third]{ITMO University, Kronverksky av., 49, 197101, Saint Petersburg, Russia (e-mail: nikona@yandex.ru)}

\begin{abstract}                
Problem of an adaptive state observer design for nonlinear system with unknown time-varying parameters and under condition of delayed measurements is considered. State observation problem was raised by many researchers (see for example \citet{b4}). In this paper the results proposed in \citet{b2}, \citet{b8}, \citet{b9}, \citet{b10} are developed. The problem is solved under assumption that the state matrix can be represented as sum of known and unknown parts. The output vector is measured with a known constant delay. An adaptive observer which reconstructs unknown state and unknown time-varying parameter is proposed.
\end{abstract}


\end{frontmatter}

\section{Introduction}
It is well known that measurements obtained with real devices are affected by delays. It makes the design of observers and control laws for dynamical systems more complicated. Many authors explored this problem recently. For the case with linear time invariant (LTI) systems, the solution of this problem is well known \citep{b1}. However, in case of linear time variant and nonlinear systems it is still widely open - see the references and literature review in \citet{b4}, \citet{b2}, \citet{b3}. There are different problems related to systems with delays considered by researchers currently. For instance, in \citet{b12} authors provide an observer-based algorithm for estimation of the unknown delay and the states of the system. A solution for the problem of state observers design for systems with constant delay was presented in \citep{b16} and state observers for the systems with time-varying delay were proposed in \citet{b2}, \citet{b8}, \citet{b9} ets. There are many examples of real systems with delayed measurements. For example, the system considered in \citet{b13, b14} is a bioreactor used for the biological treatment of wastewater. The problem of biomass regulation was solved using proportional-integral (PI) observer-based control. The slide mode controller based delayed output estimator that tracks an autonomous underwater vehicle position under time-varying disturbances and parametric variations was proposed in \citet{b15}. A discrete-time case is considered in \citet{b17}. The authors proposed a fault estimation method for linear systems with known constant time-delay and unknown disturbances and measurement noise. 

In this paper we consider a nonlinear system with unknown time-varying parameters. The unknown parameters observer is designed under condition that the output signal measurements contain a known constant delay. State observer for a nonlinear system is designed based on the obtained time-varying parameter.

\section{Problem formulation}

In this paper we consider a nonlinear system described by the following equation
\begin{align}
	\label{sys}
	\dot{x}(t)&=f(x(t), u(t))+ \sum\limits_{i=1}^n [\theta_{i}(t)Q_{i}x(t)],\\ 
	y(t)&=C(t) x(\phi(t)),
	\label{out}
\end{align}
where $x(t) \in \rea^n$ is an unknown state vector, $ f(x(t), u(t)) \in \rea^n$ and  $ C(t) \in {\rea}^{n \times n}$ are known matrices and their entries are assumed to be continuous and bouded, $Q_{i} \in {\rea}^{n \times n}$ is a known matrix, $ \phi (t) $ is a continuous known nonnegative function which defines the measurement delay
\begin{equation}
		\phi (t)=t-d,\phi(t)\ge 0, 
\end{equation}
where $ d> 0 $ is a constant delay; $ \theta_{i} (t) $ is an unknown time-varying sinusoidal function defined by the equation (see for example \citet{b21}, \citet{b22})
\begin{equation} 
	\ddot{\theta_{i} }(t)=-{\omega_{i} }^{2}\theta_{i}(t), 
	\label{thet}
\end{equation}
where $ \omega_{i} > 0 $ is a constant unknown parameter.

Solving the differential equation \eqref{thet}, we can rewrite it in the following form:
\begin{align}
	\label{theta}
	{\theta_{i} }(t)={{a }}_{1i}\sin {\omega_{i} }t+{{a }}_{2i}\cos {\omega_{i} }t,
\end{align}
where $a_{1i}$ and $a_{2i}$ are corresponding constant unknown parameters.

\begin{assumption}
We suppose that $C(t)$ is the identity matrix $C(t)={I}_{n \times n}$.
\end{assumption}

\begin{assumption}
We suppose that matrices $Q_{i}$ satisfy the following conditions:
\begin{equation}
Q_{i}Q_{j} = 0,
\end{equation}
where $j \neq i$;
\begin{equation}
Q_{i}=Q_{i}^T;
\end{equation}
\begin{equation}
Q_{i}^2 = k Q_{i},
\end{equation}
where $k \in \mathbb{Z}$.
\end{assumption}

\begin{assumption}
We suppose that
\begin{equation}
y^T Q_{i} y > 0
\end{equation}
\end{assumption}

\begin{assumption}
We suppose that function $f(x, u)$ can be represented in the following form
\begin{equation}
f(x(t),u(t)) = A(y(t), u(t))x(t) + B(y(t), u(t))u(t)
\end{equation}
where $A(y(t), u(t)) \in \rea ^{n \times n}; B(y(t), u(t)) \in \rea ^n$.
\end{assumption}

According to the Assumption 4 the system \eqref{sys}, \eqref{out} can be rewritten in the form
\begin{align}
	\dot{x}(t)=A(y(t), u(t))x(t) + &B(y(t), u(t))u(t)+ \\
	 &\sum\limits_{i=1}^n [\theta_{i}(t)Q_{i}x(t)], \nonumber 
\end{align}
\begin{equation}
	y(t)=C(t) x(\phi(t)).
\end{equation}

\section{Main Result}

\subsection{Preliminary transformations}

The problem of an observer design will be solved in three steps. At first, the unknown constant parameters $\omega_{i}$ will be estimated. Secondly, an observer for unknown continuous function $ \theta_{i} $ will be designed. At the last step, a state observer for the system \eqref{sys}, \eqref{out} will be designed.

Let us consider the system \eqref{sys}, \eqref{out} at the moment $ t – d $. The equation \eqref{sys} can be rewritten in the following way
\begin{align}
	\dot{x}_d=f_d+ \sum\limits_{i=1}^n [\eta_{i}Q_{i}x_d],
	\label{sys_delay}
\end{align}
where (by virtue of the Assumption 1)
\begin{subequations}
	\begin{align}
	{x}_d =x(\phi (t))&=y, f_d=f(x_d,u(\phi (t)), 
	\\ &\eta(t)=\theta (\phi (t)).
	\label{eta_theta} 
	\end{align}
	\end{subequations}
The new parameter $ \eta_{i}(t) $ is defined by an equation:
\begin{equation}
	\ddot{\eta_{i} }=-{\omega_{i} }^{2}\eta_{i} .
	\label{eta}
\end{equation}
Define function $\Psi_{j}$
\begin{equation}
\label{PSI}
\Psi_{j}= y^{T}Q_{j}y=x_{d}^{T}Q_{j}x_{d},
\end{equation}
where $j=\overline{1,n}$. 

For derivative of \eqref{PSI} we have:
\begin{equation}
\label{PSI_dot}
\dot{\Psi}_{j}=\alpha_{j}(x_{d}) + [\sum\limits_{i=1}^n \eta_{i}x_{d}^{T}Q_{i}^{T}]Q_{j}x_{d}+x_{d}^{T}Q_{j}[\sum\limits_{i=1}^n Q_{i}x_{d}\eta_{i}],
\end{equation}
where
\begin{equation}
\alpha_{j}(x_{d})=f_{d}^{T}Q_{j}x_{d}+x_{d}^{T}Q_{j}f_{d}
\end{equation}
is known.

According to the Assumption 2, the equation \eqref{PSI_dot} can be rewritten as following:
\begin{equation}
\dot{\Psi}_{j}=\alpha_{j}(x_{d})+2k\eta_{j}x_{d}^{T}Q_{j}x_{d}=\alpha_{j}(x_{d})+2k\eta_{j}\Psi_{j}.
\end{equation}

Let us define a new variable $\xi_{j}$ by the equation
\begin{equation}
	\xi_{j} =\ln \Psi_{j}.
\end{equation}

Its derivative can be written in the following form
\begin{equation}
	\dot{\xi}_{j}=\frac{\dot{\Psi}_{j}}{\Psi_{j}}=\frac{\alpha_{j} ({x}_d)}{\Psi_{j}}+2k\eta_{j} .
\end{equation}

Function $\eta_{j}$ can be obtained from the previous equation:
\begin{equation}
\label{eta_new}
\eta_{j}=\frac{1}{2k}(\dot{\xi}_{j}-\beta_{j}(x_{d})),
\end{equation}
where $\beta_{j}(x_{d})=\frac{\alpha_{j}(x_{d})}{\Psi_{j}}$.

\subsection{Estimation of the unknown parameter $\omega$}

The equation \eqref{eta} can be rewritten for $\eta_{j}(t)$ in operator form:
\begin{equation}
	{p}^{2}\eta_{j} =-{\omega_{j} }^{2}\eta_{j} ,
	\label{eta_dif}
\end{equation}
where $ p=d/dt $ is a differential operator.

Applying an LTI filter $ \frac{\lambda_1^3}{{(p+\lambda_1)}^{3}} $, where $\lambda_1>0$ to \eqref{eta_dif} we can obtain
\begin{equation}
	\label{filter}
	\frac{\lambda_1^3{p}^{2}}{{(p+\lambda_1)}^{3}}\eta_{j} =-{\omega_{j} }^{2}\frac{\lambda_1^3}{{(p+\lambda_1)}^{3}}\eta_{j} .
\end{equation}

Substituting $ \eta_{j}(t) $ from \eqref{eta_new} into the previous equation, we obtain:
\begin{equation}
\label{filt_theta}
	\frac{\lambda_1^3 {p}^{2}}{{(p+\lambda_1)}^{3}}(\dot{\xi}_{j}-\beta_{j}({x}_d))=-{\omega }^{2}\frac{\lambda_1^3}{{(p+\lambda_1)}^{3}}(\dot{\xi}_{j}-\beta_{j}({x}_d)).
\end{equation}

The equation \eqref{filt_theta} can be rewritten in the linear regression form
\begin{equation}
	q_j=\varphi_j v_j,
\end{equation}
where
\begin{equation}
	q_j=\frac{\lambda_1^3{p}^{3}}{{(p+\lambda_1)}^{3}}\xi_j +\frac{\lambda_1^3{p}^{2}}{{(p+\lambda_1)}^{3}}\beta_j({x}_d), \\
\end{equation}
\begin{equation}
	\varphi_j =\frac{\lambda_1^3}{{(p+\lambda_1)}^{3}}\beta_j({x}_d)-\frac{\lambda_1^3p}{{(p+\lambda_1)}^{3}}\xi_j ,\\
\end{equation}
\begin{equation}
	v_j={\omega_j }^{2}.
\end{equation}

For the unknown parameter $ v_j $ estimation standard gradient algorithm \citep{Ljung:99, Sastry:90} can be used:
\begin{equation}
	\dot{\hat{v}}_j=\gamma_1j \varphi_j (q_j-\varphi_j \hat{k}_j),
\end{equation}
where $\gamma_1j>0$ is an adaptation gain.

The unknown constant parameter $ \omega_j $ can be obtained as
\begin{equation}
	\label{hat_omega}
	\hat{\omega}_j =\sqrt{\left| \hat{v}_j \right|}.
\end{equation}

Consider the error
\begin{align}
	\tilde{v}_j=\hat{v}_j-v_j.
\end{align}

Then for derivative of $\tilde{v}_j$ we have
\begin{align}
	\nonumber \dot{\tilde{v}}_j=\dot{\hat{v}}_j=\gamma_1j \varphi_j &(q_j - \varphi_j \tilde{v}_j) =  \\
	\gamma_1j &\varphi_j^2 {v_j} -\gamma_1j \varphi_j^2 \hat{k}_j = - \gamma_1j \varphi_j^2 \tilde{v}_j.
\end{align}

The solution for $\tilde{v}_j$ takes the form
\begin{align}
	\tilde{v}_j(t) = \tilde{v}_{j0} e^{-\gamma_1j \int\limits_0^t \varphi_j^2d\tau}. 
\end{align}

Then
\begin{align}
	\hat{v}_j(t)=v+\tilde{v}_{j0} e^{-\gamma_1j \int\limits_0^t \varphi_j^2d\tau}
\end{align}
and for $\hat{\omega}_j$ we obtain
\begin{align}
\hat{\omega}_j(t)=\sqrt{	v+\tilde{v}_{j0} e^{-\gamma_1j \int\limits_0^t \varphi_j^2d\tau}}=\omega_j+\varepsilon_j(t),
\end{align}
where $\varepsilon_j(t)$ is exponentially  decaying term due to \\
$\tilde{v}_{j0} e^{-\gamma_1j \int\limits_0^t \varphi_j^2d\tau}$.

\subsection{Time-varying parameter observer}

The solution of \eqref{eta_dif} is a harmonical signal. At first let us design the observer of the signal $\eta_j$ under assumption that $\omega_j$ is known: 
\begin{equation}
	\eta_i(t)=a_{i1}\sin(\omega_j \phi(t))+a_{i2}\cos(\omega_i \phi(t)),
	\label{eta_sol}
\end{equation}
where $ a_{i1} $ and $a_{i2}$ are unknown constant parameters. Let us denote
\begin{equation}
	\chi_i=\begin{bmatrix}
		\chi_{i1}\\ 
		\chi_{i2} \end{bmatrix}
	=\begin{bmatrix}
		\sin (\omega_i \phi(t))\\ 
		\cos (\omega_i \phi(t)) \end{bmatrix}  .
	\label{chi}
\end{equation}

Now we can rewrite (\ref{eta_sol}) in the following form:

\begin{equation}
	\eta_i (t)=a_i^\top \chi_i ,
	\label{eta_matr}
\end{equation}
where $a_i =\begin{bmatrix}
	{a}_{i1}\\ 
	{a }_{i2} \end{bmatrix}.$

Substituting (\ref{eta_matr}) to (\ref{sys_delay}), we obtain:
\begin{equation}
	\dot{{x}}_d=f_d+\sum\limits_{i=1}^n {a_i}^\top \chi_i Q_i {x}_d.
\end{equation}

Let us apply a filter $ \frac{\lambda_2 }{p+\lambda_2 } $ to previous equation. After that the initial system can be transformed into linear regression form with $2n$ unknown parameters
\begin{equation}
	\caly =\sum\limits_{i=1}^n \left[a_{i1}\chi_{i1}\psi_i +a_{i2}\chi_{i2}\psi_i \right] 
\end{equation}
where
\begin{equation}
	\nonumber
	\caly = \frac{\lambda_2 p}{p+\lambda_2 }{x}_{d} -\frac{\lambda_2 }{p+\lambda_2 }f_d= \begin{bmatrix} {\caly }_{1} \\ \rotatebox{90}{$\ldots$} \\
		{\caly }_{n}\end{bmatrix},
\end{equation}
\begin{equation}
	\nonumber
	{\psi}_{i}= \begin{bmatrix} {\psi }_{i1} \\ \rotatebox{90}{$\ldots$} \\
		{\psi }_{in}\end{bmatrix}= Q_i x_d
\end{equation}

The linear regression model can be rewritten as a system of linear equations
\begin{equation}
	\begin{cases} 
		\label{lin_reg}
		\caly_{1}&=\sum\limits_{i=1}^n \left[a_{i1}\chi_{i1}\psi_{i1} +a_{i2}\chi_{i2}\psi_{i1} \right] \\ 
		&\rotatebox{90}{$\ldots$} \\
		{\caly }_{n}&=\sum\limits_{i=1}^n \left[a_{i1}\chi_{i1}\psi_{in} +a_{i2}\chi_{i2}\psi_{in} \right].
	\end{cases}
\end{equation}
In force of the Assumption 2 every equation of the system \eqref{lin_reg} can be rewritten in the following form:
\begin{equation}
\caly_i = \cala_i^T \Phi_i
\end{equation}
where:
\begin{equation}
\cala_i = a_i
\end{equation}
\begin{equation}
\Phi_i = \begin{bmatrix}
\chi_{i1} \psi_{i1} \\
\chi_{i2} \psi_{i2}
\end{bmatrix}
\end{equation}

For estimation of the unknown vector $\cala_i$ we suggest to use dynamical regression extention and mixing (DREM) technology \citep{b7, b11, b20} as it was developed in \cite{b8}. The observer can be written in the following form: 
\begin{align}
	\nonumber
		\dot Y_i = - \lambda_3 Y_i & +  \lambda_3\Phi_i^\top \caly_i,  \dot\Omega_i =- \lambda_3 \Omega_i +  \lambda\Phi_i^\top \Phi_i,	\\
		\dot{\hat{\cala_i}} &=-\gamma_2 \Delta_i (\Delta_i \hat{\cala_i}-\calz_i),
\end{align}

with $\lambda _3>0$ and $\gamma_2>0$, $\hat{\cala_i}=\hat{a}_i$ with the definitions
\begin{align}
\nonumber
		\calz_i = \adj\{\Omega_i\}Y_i, \Delta_i =\det\{\Omega_i\}.
\end{align}

Substituting $ {\hat{a }}_{i1}, {\hat{a }}_{i2} $ into  (\ref{eta_sol}), we have
\begin{equation}
	\hat{\eta }_i(t)={\hat{a }}_{i1}\sin (\omega_i \phi(t) )+{\hat{a }}_{i2}\cos (\omega_i  \phi(t)).
\end{equation}

We can obtain the unknown time-varying parameter $ \theta_i(t) $ by substitution of $ {\hat{a }}_{i1}, {\hat{a }}_{i2} $ and $\hat{\omega}_i$ (the proof is given in \citet{b10})
\begin{align}
	\hat{\theta}_i(t)&=\hat{\eta }_i(t+d)={\hat{a }}_{i1}\sin (\hat{\omega }_it)+{\hat{a }}_{i2}\cos (\hat{\omega }_i t). 
	\end{align}

\subsection{State observer for the system \eqref{sys} \eqref{out}}
Estimates of the state vector of \eqref{sys}, \eqref{out} can be found using generalized parameter estimation-based observer (GPEBO) technique \citep{b18} for systems with measurements delay from \citet{b2}. If we use the estimate of $\hat{\theta}_i(t)$ then we have the following identification algorithm.

By force of the Assumption 4 we consider the system \eqref{sys}, \eqref{out} in the form
	\begin{subequations}
		\label{gpebodyn}
		\begin{align}
			\label{dotxi}
			\dot{\xi}(t)&=A(t)\xi(t)+B(t)u(t)+ \sum\limits_{i=1}^n [\hat{\theta}_i Q_i] \xi(t),\\
			\label{dotphi}
			\dot{\Phi}_A(t)&=\left( A(t)+\sum\limits_{i=1}^n [\hat{\theta}_i Q_i] \right) \Phi_A(t),\;\Phi_A(0)=I_n,
		\end{align}
	\end{subequations}
	and the gradient parameter estimator
	\begin{align}
		\label{dothatthe}
		\dot {\hat e}(t) &=-\gamma_3 	\calp(t) [	\calp \hat{e}(t)-\calr(t)],
	\end{align}
	with $\gamma_3>0$.
	
	Let us define $\calr(t)$ and $\calp(t)$ as
	\begin{subequations}
		\label{gpebodyn1}
		\begin{align}
			\label{calz}
			\calr(t) &:= \adj\{\Phi_A(\phi(t))\}g(t),\\
			\calp(t)&:=\det\{\Phi_A(\phi(t))\},
			\label{del}
		\end{align}
	\end{subequations}
	where $\adj\{\cdot\}$ is the adjugate matrix. 
	
	Then we define the state estimate as
	\begin{subequations}
	\begin{align}
		\label{ftcgpebo}
		\hat x(t) &= \xi(t) {+} \Phi_A(t) \hat{e}_{FT},\\
		\label{e_ft}
		\hat{e}_{FT} &= {1 \over 1 - w_c(t)}[\hat e(t) - w_c(t) \hat e(0)]
	\end{align}
	\end{subequations}
	with
	\begin{align}
		\dot w(t)  &= -\gamma \calp^2(t) w(t), \; w(0)=1,
		\label{dotw}
	\end{align}
	and $w_c(t)$ defined via the clipping function
	\begin{equation}
		\label{wc}
		w_c(t) = \left\{ \begin{array}{lcl} w(t) & \;\mbox{if}\; & w(t) \leq 1-\mu, \\ 1-\mu & \;\mbox{if}\; & w(t) > 1-\mu, \end{array} \right.
	\end{equation}
	where $\mu \in (0,1)$ is a designer chosen parameter.

	Then $\hat x(t) = x(t),\;\forall t \geq t_c,$ for some $t_c \in (0,\infty)$

\begin{pf}
	
Consider error equation
\begin{align}
	\label{error}
	e(t)=x(t)-\xi(t).
\end{align}
Then, taking into account $\hat{\theta}(t)=\theta(t)$ we have
\begin{align}
\dot{e}(t)=\left( A(t)- \theta(t) I_n  \right) e(t).
\end{align}
The solution for $e(t)$ can be found in the following form
\begin{align}
	\label{e}
	e(t)=\Phi_A e(0),
\end{align}
where $e(0)=x(0) - \xi(0)$ and $\Phi_A$ is the fundamental matrix defined by equation \eqref{dotphi}. For zero initial conditions in \eqref{dotxi} we have  $e(0)=x(0)$.

After substitution \eqref{e} into \eqref{error} we can write
\begin{align}
	\label{error_1}
	\Phi_A(t) e(0)=x(t)-\xi(t).
\end{align}
In \eqref{error_1} the unknown state vector $x$ is used, but we can implement equation written in the following form
\begin{align}
	\Phi_A(\phi(t)) e(0)=x(\phi(t))-\xi(\phi(t)),
\end{align}
and we obtain linear regression equation
\begin{align}
g(t)=	\Phi_A(\phi(t)) e(0),
\end{align}
where $g(t)=x(\phi(t))-\xi(\phi(t))$. 

Now we can find vector of initial conditions $(\hat{e}(0))$ using gradient approach \eqref{dothatthe}  or $(\hat{e}_{FT}(0))$ using finite time  algorithm \eqref{e_ft}. After estimation of initial condition the state vector can be found by \eqref{ftcgpebo}.
\end{pf}

\section{Simulation Results}

Consider system \eqref{sys}, \eqref{out} with the following parameters:
\begin{equation}
\nonumber
f(x(t), u(t))= A(y(t), u(t))x(t)+B(t)u(t), 
\end{equation}
where $A(y(t), u(t))= \begin{bmatrix}
0  &0.1-0.1\sin(t)\\ 
-1   &  -1+0.5\cos(2t)
\end{bmatrix}$; \\
$B(y(t), u(t)) = \begin{bmatrix}
-1 \\ 4
\end{bmatrix}$;
$u = 2\sin(t)$;
 $C=
\begin{bmatrix}
	1 \; &0\\0 \;&1 
\end{bmatrix}$; \\
$Q_1=
\begin{bmatrix}
	1 \; &0\\0 \;&0 
\end{bmatrix}$;
$Q_2=
\begin{bmatrix}
	0 \; &0\\0 \;&1 
\end{bmatrix}$, $\omega_1=5$, $\omega_2=5$.
For simulation the following initial conditions were used $x(0)=
\begin{bmatrix}
	-5 \\ 5
\end{bmatrix}$, 
$\theta_1(0)=3, \dot{\theta}_1(0)=1.5,\theta_2(0)=2, \dot{\theta}_2(0)=1$. 
As filter parameters $\lambda_1=\lambda_2=\lambda_3=1$ were used. The constant delay in output signal is $d=2$.

\begin{figure}[hbtp]
	\begin{center}
		\includegraphics[width=1 \linewidth]{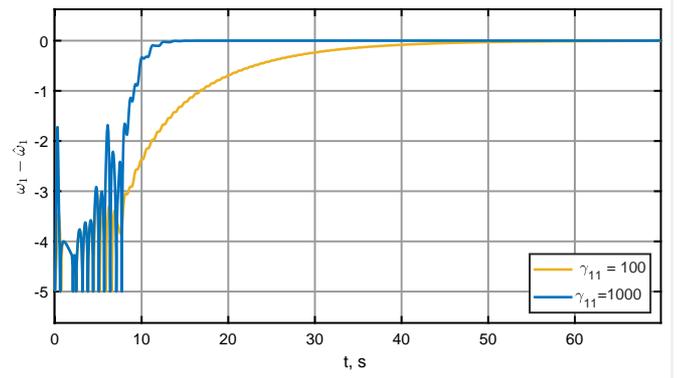}
		\caption{Transients of the error $\tilde{\omega}_1=\omega_1-\hat{\omega}_1$ for different values of adaptation gain $\gamma_{11}$}
		\label{omega1err}
	\end{center}
\end{figure}

\begin{figure}[hbtp]
	\begin{center}
		\includegraphics[width=1 \linewidth]{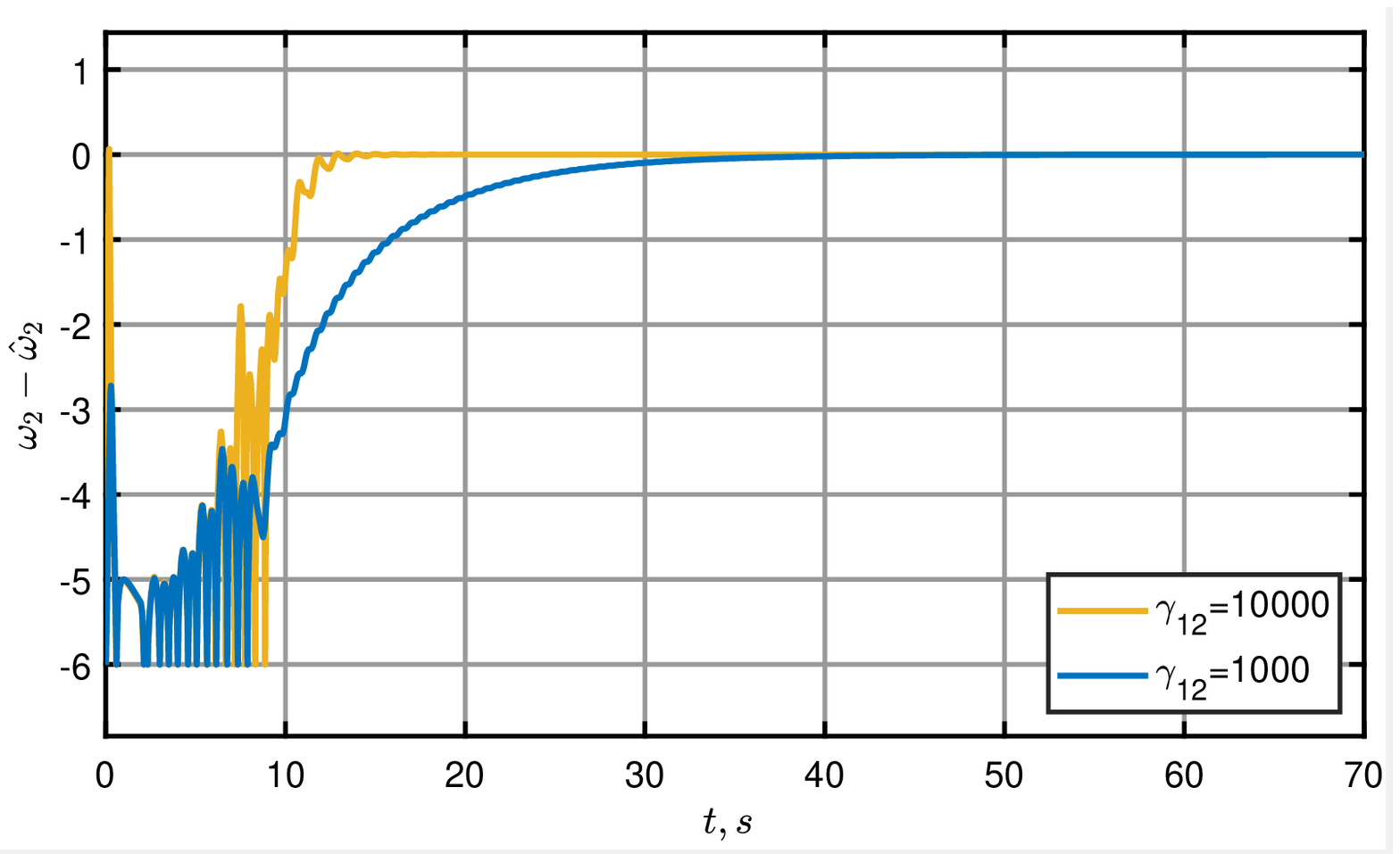}
		\caption{Transients of the error $\tilde{\omega}_2=\omega_2-\hat{\omega}_2$ for different values of adaptation gain $\gamma_{12}$}
		\label{omega2err}
	\end{center}
\end{figure}

\begin{figure}[hbtp]
	\begin{center}
		\includegraphics[width=1 \linewidth]{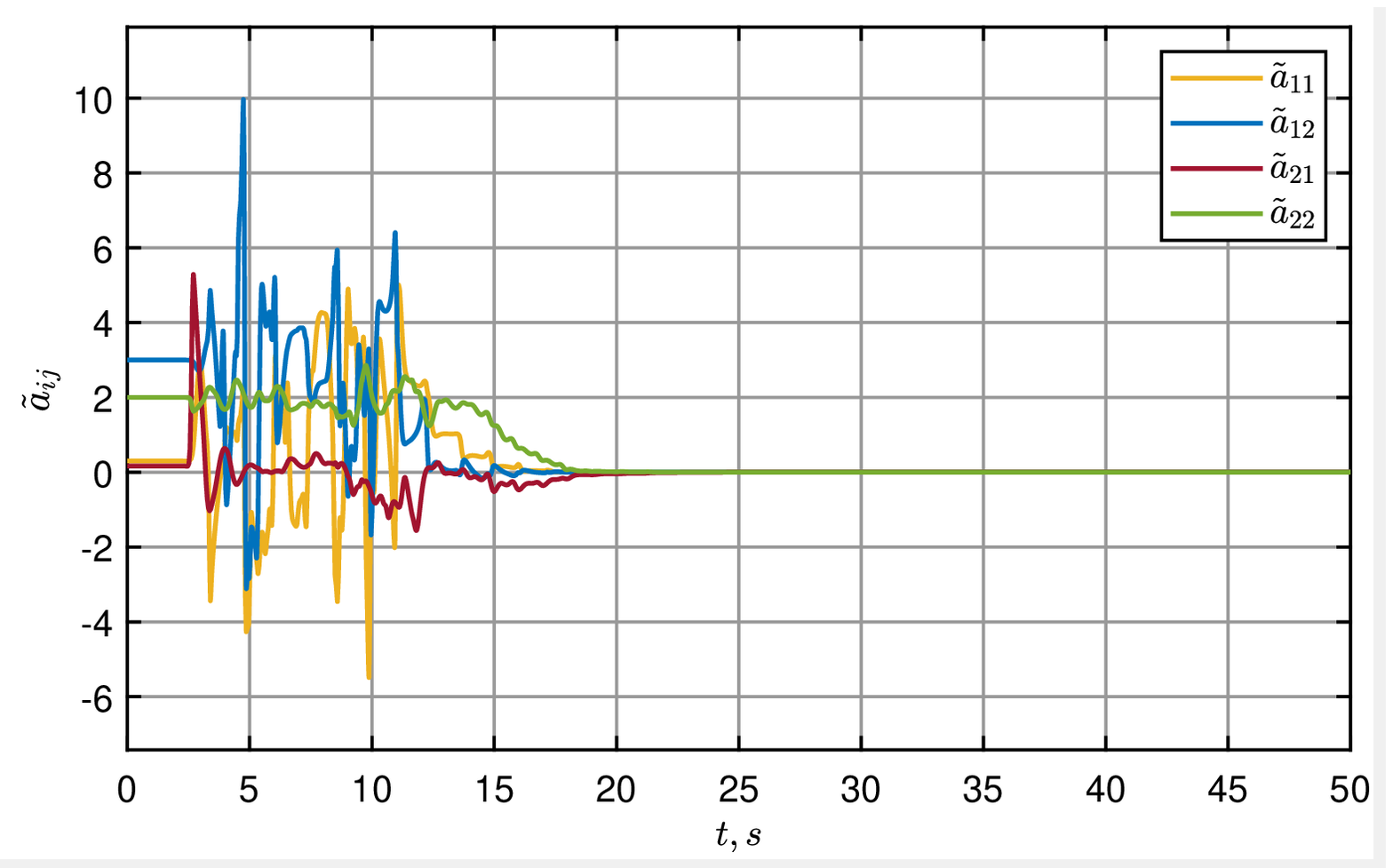}
		\caption{Transients of the error $\tilde{a}_{ij}=a_{ij}-\hat{a}_{ij}$ when the adaptation gain is $\gamma_2 = 1$}
		\label{aijerr1}
	\end{center}
\end{figure}

\begin{figure}[hbtp]
	\begin{center}
		\includegraphics[width=1 \linewidth]{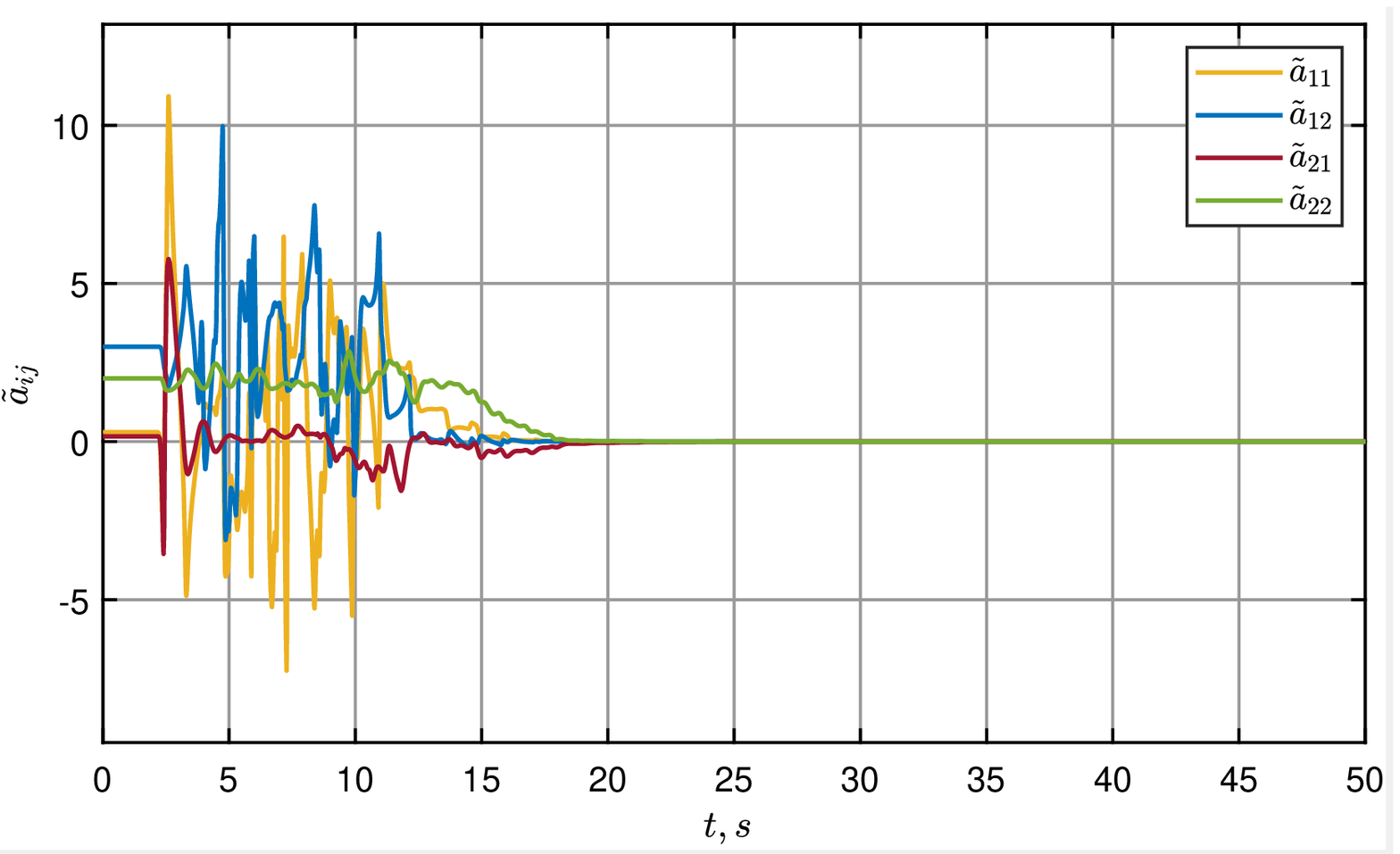}
		\caption{Transients of the error $\tilde{a}_{ij}=a_{ij}-\hat{a}_{ij}$ when the adaptation gain is $\gamma_2 = 100$}
		\label{aijerr2}
	\end{center}
\end{figure}

\begin{figure}[hbtp]
	\begin{center}
		\includegraphics[width=1 \linewidth]{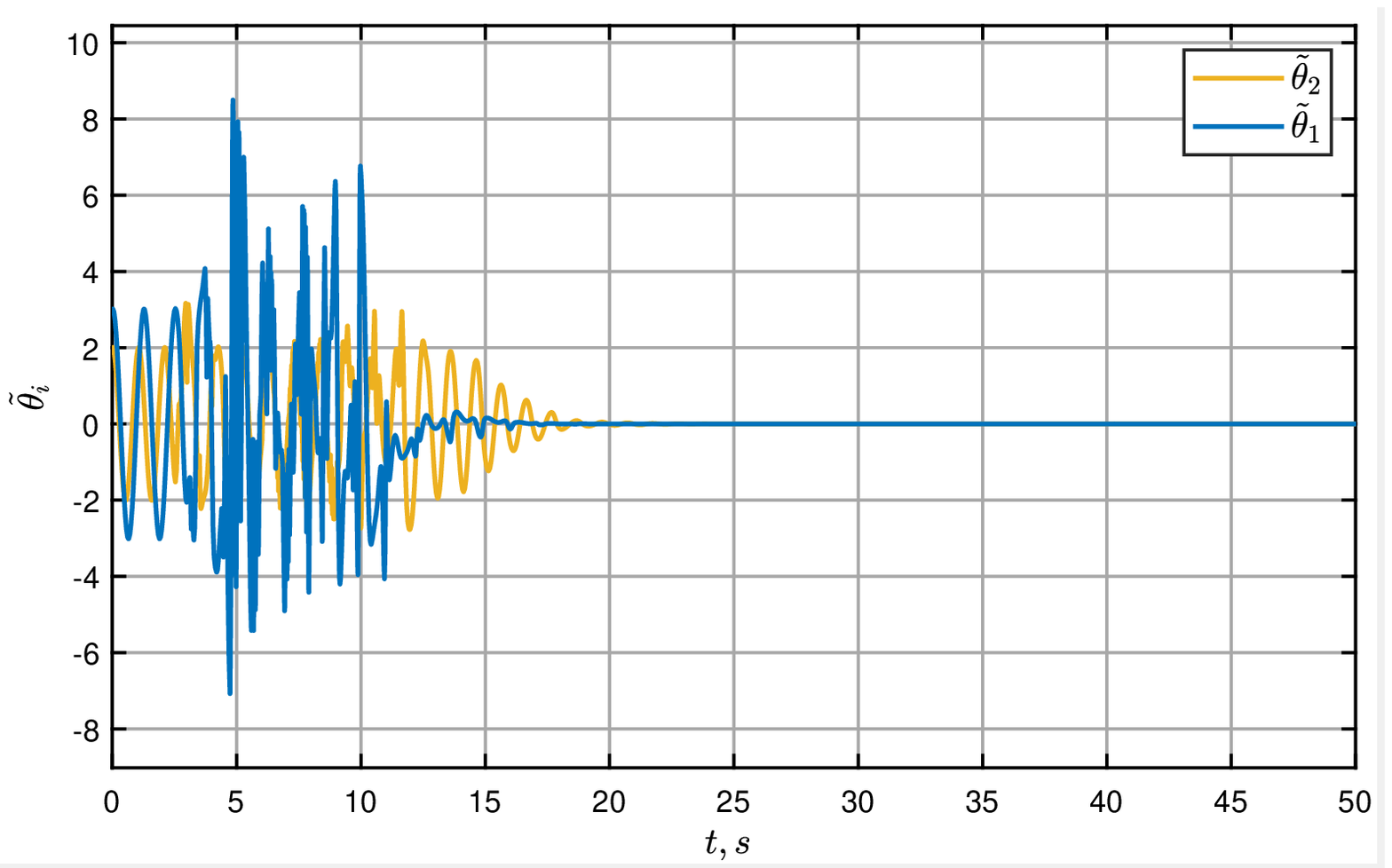}
		\caption{Transients of the error $\tilde{\theta}_{i}=\theta_{i}-\hat{\theta}_{i}$ when the adaptation gain is $\gamma_2 = 1$}
		\label{1theta12}
	\end{center}
\end{figure}

\begin{figure}[hbtp]
	\begin{center}
		\includegraphics[width=1 \linewidth]{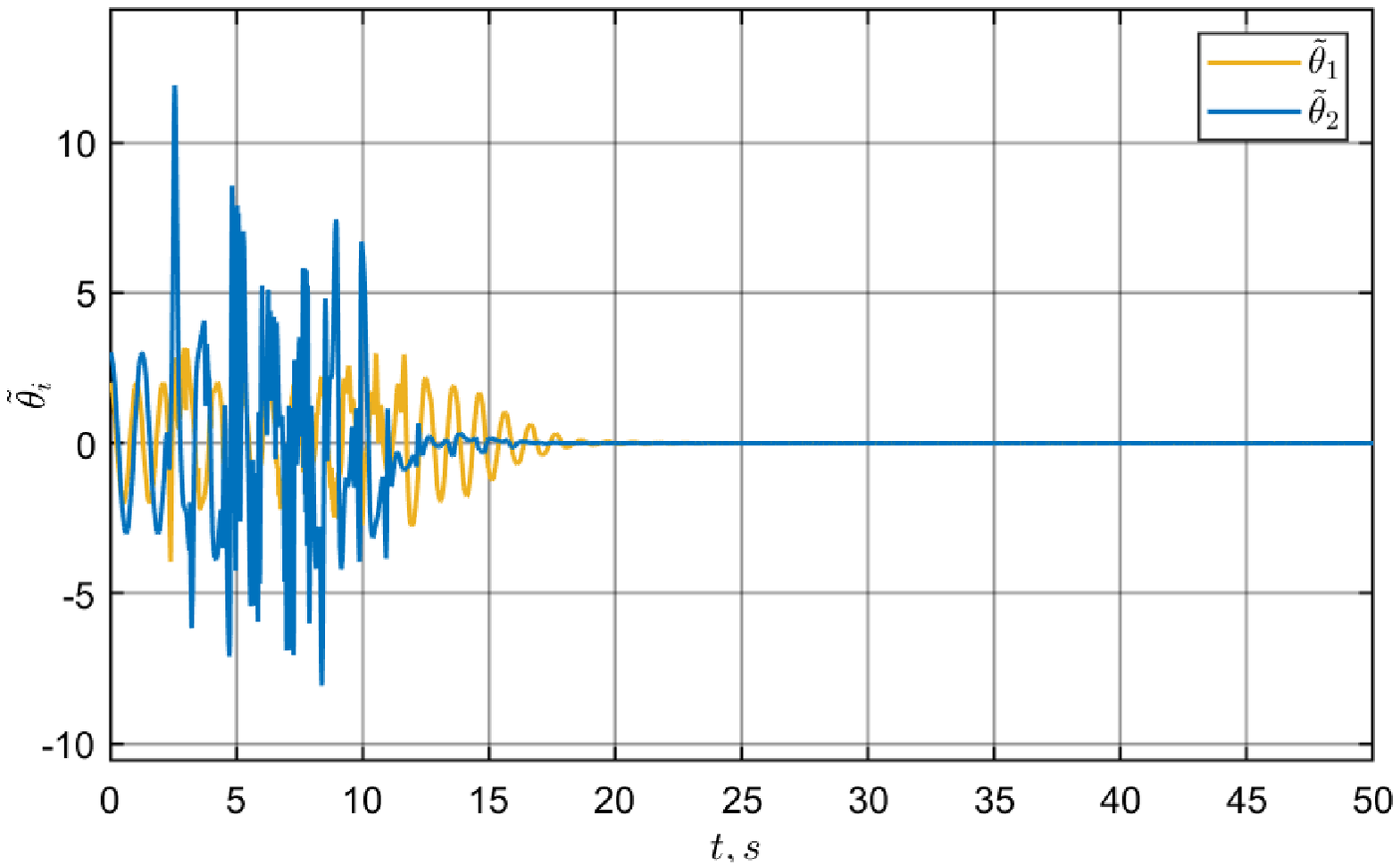}
		\caption{Transients of the error $\tilde{\theta}_{i}=\theta_{i}-\hat{\theta}_{i}$ when the adaptation gain is $\gamma_2 = 100$}
		\label{2theta12}
	\end{center}
\end{figure}

Fig.~\ref{omega1err} demonstrates transients of identification error $\tilde{\omega}_1=\theta-\hat{\omega}_1$ for different values of adaptation gain $\gamma_{11}$. Fig.~\ref{omega2err} demonstrates transients of identification error $\tilde{\omega}_2=\theta-\hat{\omega}_2$ for different values of adaptation gain $\gamma_{12}$. Fig.~\ref{aijerr1} and fig. ~\ref{aijerr2} demonstrate transients of identification errors for four unknown coefficients of unknown functions $\theta_1$ and $\theta_2$ with different values of the adaptation gain $\gamma_2$ and fixed values of adaptation gains $\gamma_{11}=1000$ and $\gamma_{12}=10000$. For the estimation of this coefficients we use the estimated values of $\omega_1$ and $\omega_2$ obtained by \eqref{hat_omega}. Fig.~\ref{1theta12} and fig.~\ref{2theta12} demonstrates transients of identification error $\tilde{\theta}_1=\theta_1-\hat{\theta}_1$ and $\tilde{\theta}_2=\theta_2-\hat{\theta}_2$ for different values of adaptation gain $\gamma_2$ and fixed values of adaptation gains $\gamma_{11}=1000$ and $\gamma_{12}=10000$. 

\begin{figure}[hbtp]
	\begin{center}
		\includegraphics[width=1 \linewidth]{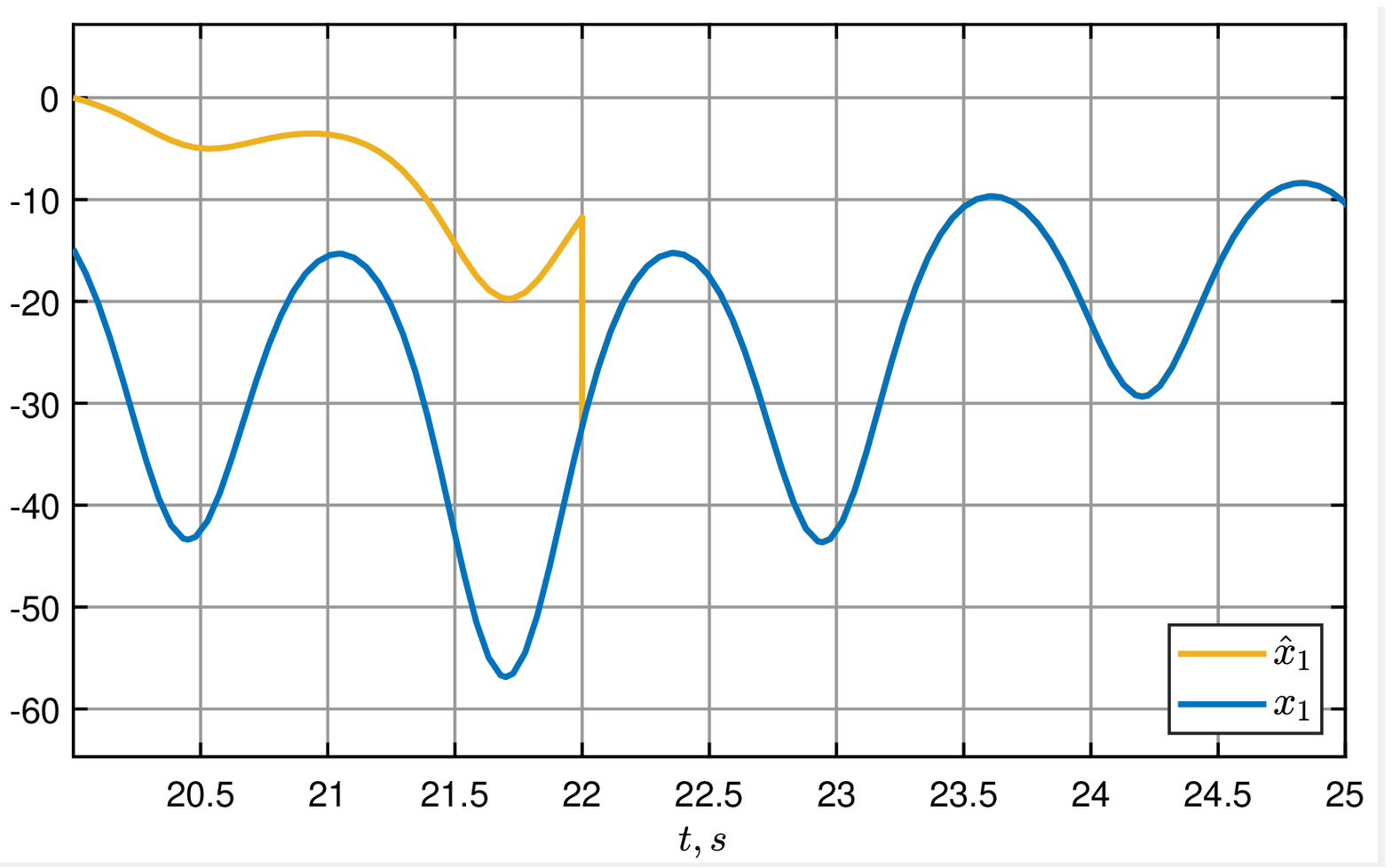}
		\caption{Transients of state variable $x_1(t)$ and its estimate $\hat x_1(t)$}
		\label{state_comp1}
	\end{center}
\end{figure}

\begin{figure}[hbtp]
	\begin{center}
		\includegraphics[width=1 \linewidth]{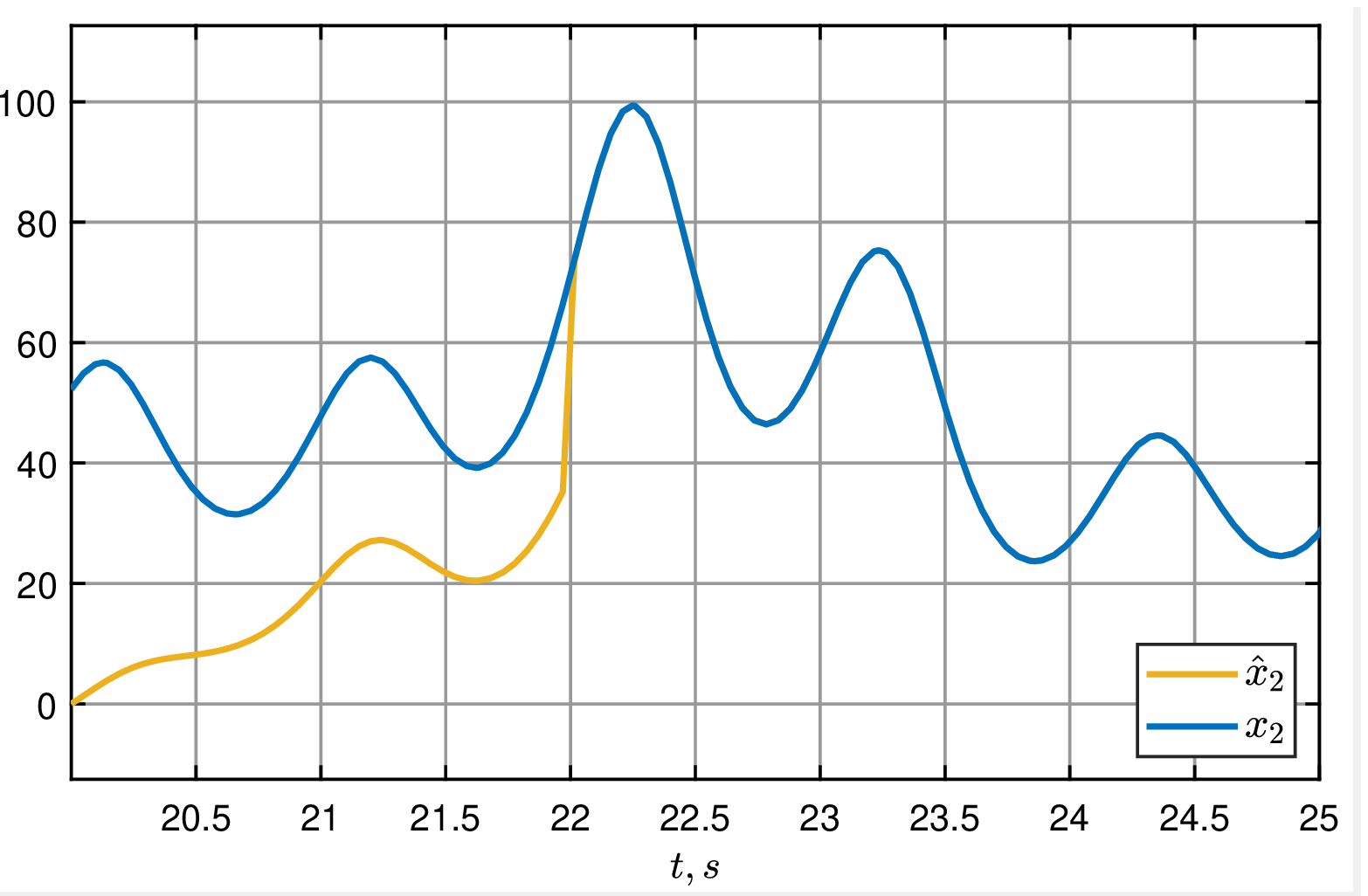}
		\caption{Transients of state variable $x_2(t)$ and its estimate $\hat x_2(t)$}
		\label{state_comp2}
	\end{center}
\end{figure}

\begin{figure}[hbtp]
	\begin{center}
		\includegraphics[width=1 \linewidth]{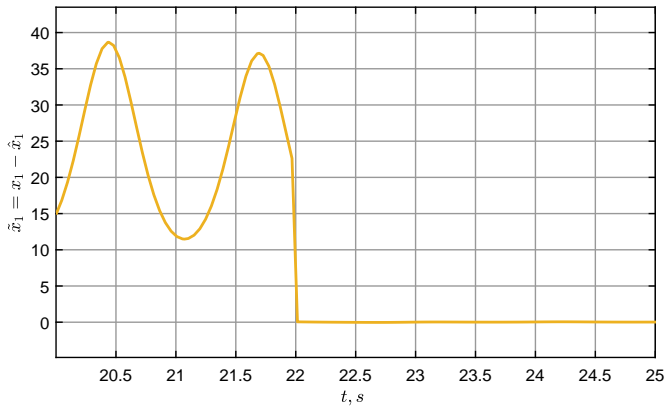}
		\caption{Transient of identification error $\tilde{x}_1=x_1-\hat{x}_1$}
		\label{state_err1}
	\end{center}
\end{figure}

\begin{figure}[hbtp]
	\begin{center}
		\includegraphics[width=1 \linewidth]{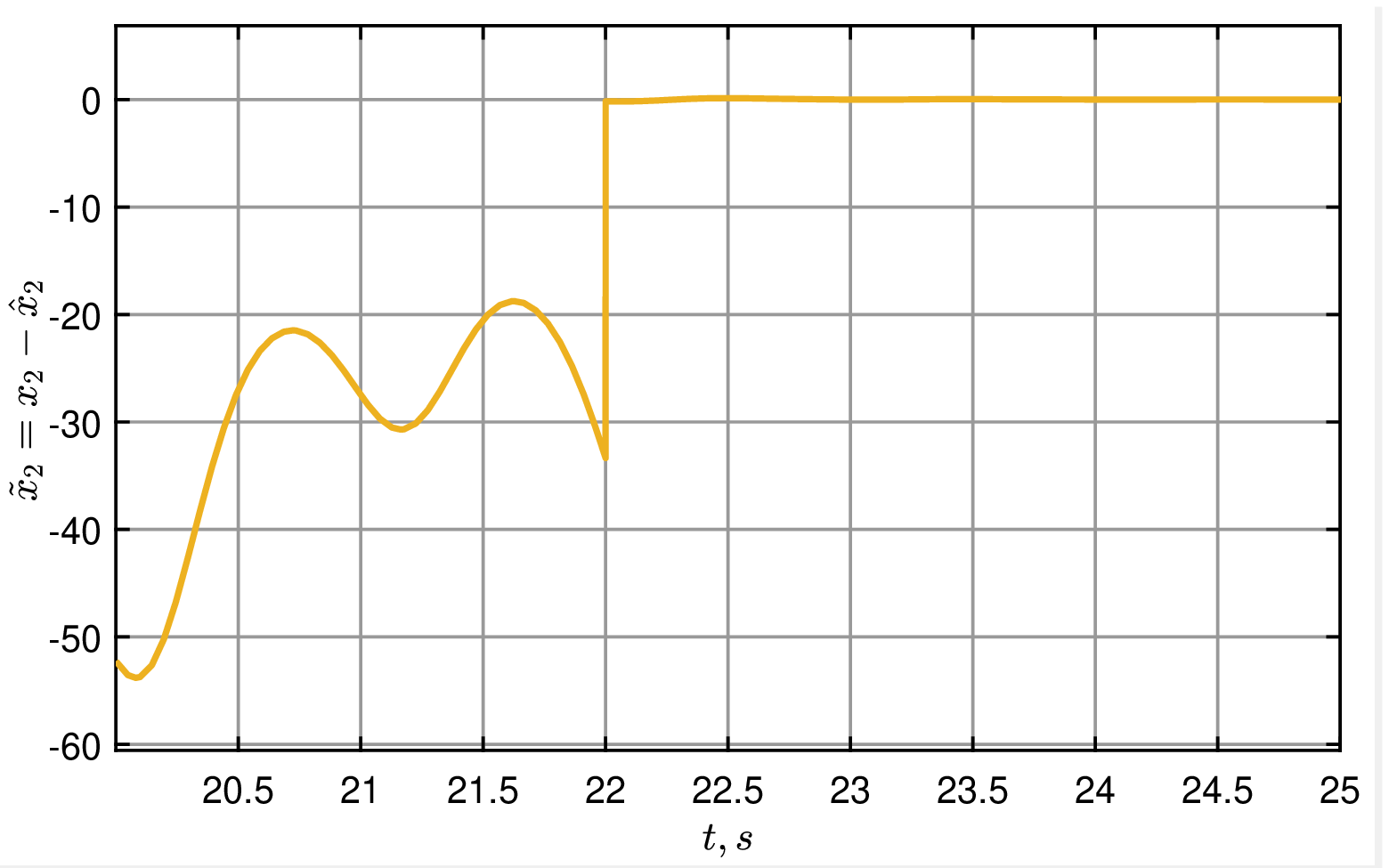}
		\caption{Transient of identification error $\tilde{x}_2=x_2-\hat{x}_2$}
		\label{state_err2}
	\end{center}
\end{figure}

As we have previously estimated the unknown time-varying parameters $\theta_1(t)$ and $\theta_2(t)$, we start estimation of the state vector. For the simulation we used $\mu=0.01$ in \eqref{wc} and the adaptation gain $\gamma_3=100$ in \eqref{dothatthe}. The state vector observer is switched on after the $\theta_1(t)$ and $\theta_2(t)$ are estimated. On the fig. \ref{state_comp1} and \ref{state_comp2} the transients of state variables $x_1(t)$ and $x_2(t)$ and their estimates $\hat x_1(t)$ and $\hat x_2(t)$ are demonstrated. On the fig. \ref{state_err1} and \ref{state_err2} the transient of identification error $\tilde{x}_1=x_1-\hat{x}_1$ and $\tilde{x}_2=x_2-\hat{x}_2$ are demonstrated. The algorithm provides convergence of the identification error to zero in finite time.


\section{Concluding Remarks} 

An adaptive observer for a nonlinear system \eqref{sys}, \eqref{out} with delayed measurements was presented. The system contains unknown time-variant parameters in the state matrix. At first we estimate constant parameters of the unknown sinusoidal function. Based on the obtained estimates we can calculate the time-varying parameters. After that we use finite time observer based on GPEBO+DREM technique to estimate unknown state vector. To demonstrate efficiency of the proposed algorithm the simulation results was presented. 

\bibliography{ifacconf}             
                                                  
\end{document}